\documentclass[twoside,reqno,12pt]{amsart}

\usepackage{latexsym}

\newcommand{\qdn}{\hspace*{-1.5mm}}
\newcommand{\qqdn}{\hspace*{-2.5mm}}
\newcommand{\xqdn}{\hspace*{-5.0mm}}
\newcommand{\xxqdn}{\hspace*{-10mm}}

\newcommand{\fns}{\footnotesize}

\newcommand{\ffnk}[4]{\left[\qdn\ba{#1}#3\\#4\ea{\!;\:#2}\right]}

\newcommand{\binm}{\binom}

\newcommand{\nnm}{\nonumber}
\newcommand{\be}{\begin{equation}}
\newcommand{\ee}{\end{equation}}
\newcommand{\ba}{\begin{array}}
\newcommand{\ea}{\end{array}}
\newcommand{\bmn}{\begin{eqnarray}}
\newcommand{\emn}{\end{eqnarray}}
\newcommand{\bnm}{\begin{eqnarray*}}
\newcommand{\enm}{\end{eqnarray*}}
\newcommand{\bln}{\begin{subequations}}
\newcommand{\eln}{\end{subequations}}

\newcommand{\lam}{\lambda}

\newtheorem{thm}{Theorem}
\newtheorem{lemm}[thm]{Lemma}

\newtheorem{entry}{Entry}

\newcommand{\bbtm}[4]{\bibitem{kn:#1}{#2,}~{#3,}~{#4.}}
\newcommand{\cito}[1]{\cite{kn:#1}}
\newcommand{\citu}[2]{\cite[#2]{kn:#1}}

\usepackage{indentfirst}

\begin{document} 
{\fns 
\title{Several transformation formulas for basic hypergeometric series}
\author{$^1$Chuanan Wei, $^2$Dianxuan Gong$^{*}$}
\dedicatory{
$^1$School of Biomedical Information and Engineering\\
  Hainan Medical University, Haikou, China\\
  $^2$College of Sciences\\
  North China University of Science and Technology, Tangshan, China\\
   }
\thanks{$^{*}$Corresponding author. \emph{Email addresses}:
    weichuanan78@163.com (C. Wei), gdx216@163.com (D. Gong)}

\address{ }
\footnote{\emph{2010 Mathematics Subject Classification}: Primary
05A19 and Secondary 33D15.}

 \keywords{Basic hypergeometric series;
 Ramanujan's ${_1\psi_1}$ summation formula;
Andrews' identity;
   Ramanujan's reciprocity theorem}

\begin{abstract}
In 1981, Andrews gave a four-variable generalization of Ramanujan's
${_1\psi_1}$ summation formula. We establish a six-variable
generalization of Andrews' identity according to the transformation
formula for two ${_8\phi_7}$ series and Bailey's transformation
formula for three ${_8\phi_7}$ series. Then it is used to find a
six-variable generalization of Ramanujan's reciprocity theorem,
which is different from Liu's formula. We derive the generalizations
of Bailey's two $_3\psi_3$ summation formulas in terms of two
limiting relations and Bailey's another transformation formula for
three $_8\phi_7$ series. Based on the two limiting relations, some
different results involving bilateral basic hypergeometric series
are also deduced from the Guo--Schlosser transformation formula and
other two transformation formulas.

\end{abstract}

\maketitle\thispagestyle{empty}
\markboth{Chuanan Wei, Dianxuan Gong}
         {Andrews' identity and Ramanujan's reciprocity theorem}

\section{Introduction}
\hspace{0.3cm} Let $x$, $q$ be complex numbers  with $|q|<1$ and $n$
a nonnegative integer. Define the $q$-shifted factorial to be
 \bnm
(x;q)_{\infty}=\prod_{i=0}^{\infty}(1-xq^i),\quad
(x;q)_n=\frac{(x;q)_{\infty}}{(xq^n;q)_{\infty}}.
 \enm
For shortening many of the formulas in this paper, we adopt the
notation
 \bnm
(x_1,x_2,\dots,x_r;q)_{m}=(x_1;q)_{m}(x_2;q)_{m}\cdots(x_r;q)_{m},
 \enm
where $m\in\mathbb{Z}\cup\{\infty\}$. Following Gasper and Rahman
\cito{gasper}, define the bilateral basic hypergeometric series by
 \bnm\quad
{_{r}\psi_s}\ffnk{cccccc}{q,z}{a_1,a_2,\ldots,a_r}{b_1,b_2,\ldots,b_s}
  =\sum_{k=-\infty}^{\infty}\frac{(a_1,a_2,\ldots,a_r;q)_k}{(b_1,b_2,\ldots,b_s;q)_k}
\Big\{(-1)^kq^{\binm{k}{2}}\Big\}^{s-r}z^k.
 \enm
Then Ramanujan's ${_1\psi_1}$ summation formula (cf.
\citu{gasper}{Appendix (II.29)}) reads
 \bmn\label{ramanujan}
{_1\psi_1}\ffnk{cccccc}{q,z}{a}{b}=
\frac{(q,b/a,az,q/az;q)_{\infty}}{(b,q/a,z,b/az;q)_{\infty}},
 \emn
provided $|b/a|<|z|<1$. Bailey's two $_3\psi_3$ summation formulas
 (cf. \citu{gasper}{Exercise 5.18}) can be expressed as
 \bmn\label{p33-a}
\xqdn\qqdn{_3\psi_3}\ffnk{cccc}{q,\frac{q}{bcd}}{b,c,d}
 {q/b,q/c,q/d}=\frac{(q,q/bc,q/bd,q/cd;q)_{\infty}}{(q/b,q/c,q/d,q/bcd;q)_{\infty}},
 \emn
where $|q/bcd|<1$,
 \bmn \label{p33-b}
\quad\:{_3\psi_3}\ffnk{cccc}{q,\frac{q^2}{bcd}}{b,c,d}
 {q^2/b,q^2/c,q^2/d}=\frac{(q,q^2/bc,q^2/bd,q^2/cd;q)_{\infty}}{(q^2/b,q^2/c,q^2/d,q^2/bcd;q)_{\infty}},
 \emn
provided $|q^2/bcd|<1$. Some related results can be seen in
\cito{zhang}.

\hspace{0.3cm} When $b_s=q$, the bilateral basic hypergeometric
series reduces to the unilateral basic hypergeometric series
 \bnm\qquad
{_{r}\phi_{s-1}}\ffnk{cccccc}{q,z}{a_1,a_2,\ldots,a_r}{b_1,\ldots,b_{s-1}}
  =\sum_{k=0}^{\infty}\frac{(a_1,a_2,\ldots,a_r;q)_k}{(q,b_1,\ldots,b_{s-1};q)_k}
\Big\{(-1)^kq^{\binm{k}{2}}\Big\}^{s-r}z^k.
 \enm
Habitually, the very-well-poised series
 \bnm
{_{r+1}\phi_r}\ffnk{cccccc}{q,z}{a_1,q\sqrt{a_1},-q\sqrt{a_1},a_4,\ldots,a_{r+1}}{\sqrt{a_1},-\sqrt{a_1},a_1q/a_4,\ldots,a_1q/a_{r+1}}
 \enm
is usually denoted by the symbol
 \bnm
 _{r+1}W_{r}(a_1;a_4,\ldots,a_{r+1};q,z).
  \enm
Thus the transformation formula for two ${_8\phi_7}$ series (cf.
\citu{gasper}{Appendix (III.23)}) and Bailey's transformation
formula for three ${_8\phi_7}$ series (cf. \citu{gasper}{Appendix
(III.37)}) read
 \bmn\label{two-term}
  &&{_8W_7}(a;b,c,d,e,f;q,a^2q^2/bcdef)
 \nnm\\[1mm]
 &&\:=
\frac{(aq,aq/ef,\lam q/e,\lam q/f;q)_{\infty}}
 {(aq/e,aq/f,\lam q,\lam q/ef;q)_{\infty}}
 {_8W_7}(\lam;\lam b/a,\lam c/a,\lam d/a,e,f;q,aq/ef)
 \emn
with $\lam=qa^2/bcd$ and $\max\{|a^2q^2/bcdef|,|aq/ef|\}<1$,
 \bmn\label{three-term}
  &&\xxqdn{_8W_7}(a;b,c,d,e,f;q,a^2q^2/bcdef)
  \nnm\\[1mm]\nnm
 &&\xxqdn\:=\frac{(aq,aq/de,aq/df,aq/ef,eq/c,fq/c,b/a,bef/a;q)_{\infty}}{(aq/d,aq/e,aq/f,aq/def,q/c,efq/c,be/a,bf/a;q)_{\infty}}\\[1mm]
 &&\xxqdn\:\times\:{_8W_7}(ef/c;aq/bc,aq/cd,ef/a,e,f;q,bd/a)
 \nnm\\[1mm]\nnm
  &&\xxqdn\:+\:\frac{b}{a}\frac{(aq,bq/a,bq/c,bq/d,bq/e,bq/f,d,e,f;q)_{\infty}}
  {(aq/b,aq/c,aq/d,aq/e,aq/f,bd/a,be/a,bf/a,def/a;q)_{\infty}}\\[1mm]
  &&\xxqdn\:\times\:\frac{(aq/bc,bdef/a^2,a^2q/bdef;q)_{\infty}}{(aq/def,q/c,b^2q/a;q)_{\infty}}
  \nnm\\[1mm]
 &&\xxqdn\:\times\:
{_8W_7}(b^2/a;b,bc/a,bd/a,be/a,bf/a;q,a^2q^2/bcdef),
 \emn
where $\max\{|a^2q^2/bcdef|,|bd/a|\}<1$. Bailey's another
transformation formula for three $_8\phi_7$ series (cf.
\citu{gasper}{Exercise 2.15}) can be stated as
 \bmn\label{bailey-three-term}
  &&\xqdn\frac{1}{a}\frac{(qa/d,qa/e,qa/f,q/ad,q/ae,q/af;q)_{\infty}}{(qa^2,ab,ac,b/a,c/a;q)_{\infty}} {_8W_7}(a^2;ab,ac,ad,ae,af;q,q^2/abcdef)
  \nnm\\[1mm]\nnm
  &&\xqdn+\:\frac{1}{b}\frac{(qb/d,qb/e,qb/f,q/bd,q/be,q/bf;q)_{\infty}}{(qb^2,ba,bc,a/b,c/b;q)_{\infty}}{_8W_7}(b^2;ba,bc,bd,be,bf;q,q^2/abcdef)\\[1mm]
  &&\xqdn+\:\frac{1}{c}\frac{(qc/d,qc/e,qc/f,q/cd,q/ce,q/cf;q)_{\infty}}{(qc^2,ca,cb,a/c,b/c;q)_{\infty}}{_8W_7}(c^2;ca,cb,cd,ce,cf;q,q^2/abcdef)
  \nnm\\[1mm]
  &&\xqdn=\:0,
 \emn
  provided $|q^2/abcdef|<1$. For $\lambda=qa^2/bcd$, there hold the two
transformation formulas (cf. \citu{gasper}{Exercise 2.13}):
 \bmn
  &&{_8W_7}\Big(\lambda;a^{\frac{1}{2}},-a^{\frac{1}{2}},\lambda b/a,\lambda c/a,\lambda d/a;q,-q\Big)
  \nnm\\[1mm]\label{equation-a}
  &&\:=\frac{(\lambda q,-\lambda q/a,qa^{\frac{1}{2}},-qa^{\frac{1}{2}};q)_{\infty}}{(aq,-q,\lambda qa^{-\frac{1}{2}},-\lambda qa^{-\frac{1}{2}};q)_{\infty}}
  {_4\phi_3}\ffnk{ccccccc}{q,-\frac{\lambda q}{a}}{a,b,c,d}{aq/b,aq/c,aq/d},
  \emn
where $|\lambda q/a|<1$,
  \bmn
  &&\,\xxqdn{_{10}W_9}\Big(\lambda;a^{\frac{1}{2}},-a^{\frac{1}{2}},(aq)^{\frac{1}{2}},-(aq)^{\frac{1}{2}},\lambda b/a,\lambda c/a,\lambda d/a;q,\lambda q/a\Big)
  \nnm\\[1mm]\label{equation-b}
  &&\,\xxqdn\:=\frac{(\lambda q,\lambda^2q/a^2;q)_{\infty}}{(\lambda q/a,\lambda^2q/a;q)_{\infty}}
  {_4\phi_3}\ffnk{ccccccc}{q,\frac{\lambda^2q}{a^2}}{a,b,c,d}{aq/b,aq/c,aq/d},
 \emn
provided $\max\{|\lambda q/a|,|\lambda^2 q/a^2|\}<1$. Recently, Guo
and Schlosser \cito{guo-b} found some new $q$-supercongruences and
the following new transformation formula:
 \bmn\label{Guo-Schlosser}
  &&{_{12}W_{11}}\Big(a;b,c,d,ab/c,ab/d,(aq/b)^{\frac{1}{2}},-(aq/b)^{\frac{1}{2}},q(a/b)^{\frac{1}{2}},-q(a/b)^{\frac{1}{2}};q,q/b\Big)
  \nnm\\[1mm]
  &&\:=\frac{(aq,ab/c,ab/d,aq/cd;q)_{\infty}}{(ab,aq/c,aq/d,ab/cd;q)_{\infty}}
  {_4\phi_3}\ffnk{ccccccc}{q,\frac{q^2}{b^2}}{b,c,d,cd/a}{cq/b,dq/b,cdq/ab}
  \nnm\\[1mm]
  &&\:+\:\frac{(aq,c,d,cdq/ab^2;q)_{\infty}}{(ab,cq/b,dq/b,cd/ab;q)_{\infty}}
  {_4\phi_3}\ffnk{ccccccc}{q,\frac{q^2}{b^2}}{b,ab/c,ab/d,ab^2/cd}{aq/c,aq/d,abq/cd},
 \emn
where $|q/b|<1$. More results related to $q$-supercongruences can be
seen in \cite{kn:guo-a,kn:guo-c}.

\hspace{0.3cm} In 1981, Andrews \cito{andrews} gave the
four-variable generalization of \eqref{ramanujan}:
 \bmn\label{andrews}
&&\frac{1}{a}\sum_{k=0}^{\infty}\frac{(-q/a,AB/ab;q)_k}{(-A/a,-B/a;q)_{k+1}}(-b)^k
-\frac{1}{b}\sum_{k=0}^{\infty}\frac{(A,-aq/B;q)_k}{(-A/b,-a;q)_{k+1}}\bigg(-\frac{B}{b}\bigg)^k
 \nnm\\[1mm]
&&\:=\:\bigg(\frac{1}{a}-\frac{1}{b}\bigg)\frac{(q,aq/b,bq/a,A,B,AB/ab;q)_{\infty}}{(-a,-b,-A/a,-A/b,-B/a,-B/b;q)_{\infty}},
 \emn
provided $|B|<|b|<1$. By applying the $q$-exponential operator to
\eqref{ramanujan}, Liu \cito{liu-a} established the reciprocity
formula:
 \bmn\label{liu-a}
&&v\sum_{k=0}^{\infty}\frac{(q/u,acuv;q)_k}{(av,cv;q)_{k+1}}v^k
-u\sum_{k=0}^{\infty}\frac{(q/v,acuv;q)_k}{(au,cu;q)_{k+1}}u^k
 \nnm\\[1mm]
&&\:=\:\frac{(v-u)(q,uq/v,vq/u,auv,cuv,acuv;q)_{\infty}}{(u,v,au,av,cu,cv;q)_{\infty}},
 \emn
where $\max\{|u|,|v|\}<1$. In the light of
$q$-Kummer--Thomae--Whipple formula (cf. \citu{gasper}{Appendix
(III.9)}):
 \bnm\quad
{_3\phi_2}\ffnk{ccccccc}{q;\frac{de}{abc}}{a,b,c}{d,e}=\frac{(e/a,de/bc;q)_{\infty}}{(e,de/abc;q)_{\infty}}
{_3\phi_2}\ffnk{ccccccc}{q;\frac{e}{a}}{a,d/b,d/c}{d,de/bc},
 \enm
provided $\max\{|de/abc|,|e/a|\}<1$, it is not difficult to
understand the equivalence of \eqref{andrews} and \eqref{liu-a}.
 Replace $u$ by $-u/a$ and $v$ by $-v/a$ in
 \eqref{liu-a} to obtain
 \bnm
&&v\sum_{k=0}^{\infty}\frac{(-aq/u,cuv/a;q)_k}{(-v,-cv/a;q)_{k+1}}\bigg(-\frac{v}{a}\bigg)^k
-u\sum_{k=0}^{\infty}\frac{(-aq/v,cuv/a;q)_k}{(-u,-cu/a;q)_{k+1}}\bigg(-\frac{u}{a}\bigg)^k
 \nnm\\[1mm]
&&\:=\:\frac{(v-u)(q,uq/v,vq/u,uv/a,cuv/a,cuv/a^2;q)_{\infty}}{(-u,-v,-u/a,-v/a,-cu/a,-cv/a;q)_{\infty}}.
 \enm
When $a\to \infty$, the last equation reduces to Ramanujan's
reciprocity theorem (cf. \citu{ramanujan}{p. 40}):
 \bmn \label{Ramanujan-a}
&&\xqdn\Big(1+\frac{1}{u}\Big)\sum_{k=0}^{\infty}\frac{q^{\binm{1+k}{2}}}{(-qv;q)_k}\bigg(-\frac{v}{u}\bigg)^k
-\Big(1+\frac{1}{v}\Big)\sum_{k=0}^{\infty}\frac{q^{\binm{1+k}{2}}}{(-qu;q)_k}\bigg(-\frac{u}{v}\bigg)^k
\nnm\\[1mm]
&&\xqdn\:=\Big(\frac{1}{u}-\frac{1}{v}\Big)\frac{(q,uq/v,vq/u;q)_{\infty}}{(-uq,-vq;q)_{\infty}}.
 \emn
The first published proof of \eqref{Ramanujan-a} was furnished by
Andrews \cito{andrews}. For three different proofs of it, the reader
is referred to Berndt, Chan, Yeap, and  Yee \cito{berndt}. By means
of a limiting case of Watson's transformation formula (cf.
\citu{gasper}{Exercise
 2.22}), Kang \cito{kang} derive the equivalent form
  of \eqref{andrews}:
\bnm
 &&\xxqdn\Big(1+\frac{1}{b}\Big)\sum_{k=0}^{\infty}
 \frac{(c,d,cd/ab;q)_k(1+q^{2k}cd/b)}{(-aq;q)_k(-c/b,-d/b;q)_{k+1}}q^{\binm{1+k}{2}}
\Big(-\frac{a}{b}\Big)^k
 \nnm\\[1mm]\nnm
&&\xxqdn\:-\:\Big(1+\frac{1}{a}\Big)\sum_{k=0}^{\infty}
 \frac{(c,d,cd/ab;q)_k(1+q^{2k}cd/a)}{(-bq;q)_k(-c/a,-d/a;q)_{k+1}}q^{\binm{1+k}{2}}
\Big(-\frac{b}{a}\Big)^k\\[1mm]
 &&\xxqdn\:=\:\Big(\frac{1}{b}-\frac{1}{a}\Big)
\frac{(q,aq/b,bq/a,c,d,cd/ab;q)_{\infty}}{(-aq,-bq,-c/a,-c/b,-d/a,-d/b;q)_{\infty}},
 \enm
which is also a four-variable generalization of \eqref{Ramanujan-a}.
By splitting a bilateral series into two unilateral series, Ma
\cito{ma} got the five-variable generalization of
\eqref{Ramanujan-a}:
 \bnm
&&\sum_{k=0}^{\infty}\Big(1-\frac{aq^{2k+1}}{b}\Big)\frac{(-1/b;q)_{k+1}}{(-aq;q)_k}
  \frac{(-aq/c,-aq/d,-aq/e;q)_k}{(-c/b,-d/b,-e/b;q)_{k+1}}\bigg(\frac{cde}{abq}\bigg)^k\\[1mm]
&&\:-\:\sum_{k=0}^{\infty}\Big(1-\frac{bq^{2k+1}}{a}\Big)\frac{(-1/a;q)_{k+1}}{(-bq;q)_k}
  \frac{(-bq/c,-bq/d,-bq/e;q)_k}{(-c/a,-d/a,-e/a;q)_{k+1}}\bigg(\frac{cde}{abq}\bigg)^k\\[1mm]
&&\:=\Big(\frac{1}{b}-\frac{1}{a}\Big)
 \frac{(q,aq/b,bq/a,c,d,e,cd/ab,ce/ab,de/ab;q)_{\infty}}{(-aq,-bq,-c/a,-c/b,-d/a,-d/b,-e/a,-e/b,cde/abq;q)_{\infty}},
 \enm
where $|cde/abq|<1$. In accordance with the method of $q$-partial
differential equations, Liu \cito{liu-b} gained recently the
six-variable generalization of \eqref{Ramanujan-a}:
 \bmn\label{liu-b}
&&\xqdn\qqdn
v\sum_{k=0}^{\infty}\frac{(q/u,acuv,bcuv;q)_k}{(av,bv,cv;q)_{k+1}}v^k
{_3\phi_2}\ffnk{ccccccc}{q;\frac{abcruv}{q}}{q^{k+1},vq^{k+1}/r,q/cu}{avq^{k+1},bvq^{k+1}}
 \nnm\\[1mm]\nnm
&&\xqdn\qqdn\:-\:u\sum_{k=0}^{\infty}\frac{(q/v,acuv,bcuv;q)_k}{(au,bu,cu;q)_{k+1}}u^k
{_3\phi_2}\ffnk{ccccccc}{q;\frac{abcruv}{q}}{q^{k+1},uq^{k+1}/r,q/cv}{auq^{k+1},buq^{k+1}}
\\[1mm]
&&\xqdn\qqdn\:=\frac{(v-u)(q,uq/v,vq/u,auv,buv,cuv,abr,acuv,bcuv;q)_{\infty}}{(u,v,au,av,bu,bv,cu,cv,abcruv/q;q)_{\infty}}
\nnm\\[1mm]
&&\xqdn\qqdn\:\times\:{_3\phi_2}\ffnk{ccccccc}{q;abr}{u,v,uv/r}{auv,buv},
 \emn
provided
$\max\{|u|,|v|,|au|,|av|,|bu|,|bv|,|cu|,|cv|,|abr|,|abcruv/q|\}<1$.

\hspace{0.3cm} Inspired by the works just mentioned, we shall
construct a six-variable generalization of Andrews' identity
\eqref{andrews} via the utilization of \eqref{two-term} and
\eqref{three-term} in Section 2. Then it is used to offer a
six-variable generalization of Ramanujan's reciprocity theorem
\eqref{Ramanujan-a}, which is essentially different from
\eqref{liu-b}, in Section 3. We shall establish two limiting
relations and utilize them to deduce some transformation formulas
involving bilateral basic hypergeometric series in Section 4.
\section{A six-variable generalization of Andrews' identity}
\hspace{0.3cm} Firstly, we shall display a six-variable
generalization of Andrews' identity \eqref{andrews} in the following
theorem.

\begin{thm}  \label{thm-a}
Let $A$, $B$, $a$, $b$, $x$, $y$ be complex numbers. Then
 \bnm
&&\xqdn\frac{1}{a}\sum_{k=0}^{\infty}\frac{(-q/a,AB/ab;q)_k(-b)^k}{(-A/a,-B/a;q)_{k+1}}
\frac{(x,y,-xa/A,-xa/B;q)_k(1-q^{2k}x)}{(q,-xa,xq/y;q)_{k}(xab/AB;q)_{k+1}}\bigg(\frac{q}{y}\bigg)^k\\[1mm]
&&\xqdn\:-\:\frac{\lambda(y)}{b}\sum_{k=0}^{\infty}\frac{(A,-ya/B;q)_k(-B/b)^k}{(-ya/q,-yA/bq;q)_{k+1}}
\frac{(xya^2/Bq,y,-xa^2b/AB,-xa/B;q)_k}{(q,-xa,xa^2/B;q)_{k}}\\[1mm]
&&\xqdn\:\times\:\frac{(1-q^{2k-1}xya^2/B)}{(xya^2/ABq;q)_{k+1}}\\[1mm]
&&\xqdn\:=\frac{\theta(y)(y,ya/b,bq/ya,A,Bq/y,AB/ab,-xa/B,-xaq/yA,-xa^2bq/yAB;q)_{\infty}}{(-a,-bq/y,-A/a,-yA/b,-B/a,-B/b,-xa,xab/AB,xa^2/AB;q)_{\infty}}\\[1mm]
&&\xqdn\:\times\:{_8W_7}(xa^2/yA;q/y,-xa/y,xa^2/AB,-a,-B/b;q,-bq/A),
 \enm
where $\max\{|bq/y|,|B/b|,|bq/A|\}<1$ and the notations
$\lambda(y)$, $\theta(y)$ given by
 \bnm
  &&\xqdn\lambda(y)=\frac{y}{q}\frac{(x,xa^2/B,-aq/A,-bq/A,-ya/q,-Bq/ya,xabq/yAB,xya^2/ABq;q)_{\infty}}
  {(xq/y,xya^2/Bq,-ya/A,-bq^2/yA,-a,-B/a,xab/AB,xa^2/AB;q)_{\infty}},\\[1mm]
   &&\xqdn\theta(y)=\frac{-1}{A+b}\frac{(x,xa^2/A,-aq/A,-b/A;q)_{\infty}}
  {(xq/y,xa^2q/yA,-ya/A,-bq/yA;q)_{\infty}}.
  \enm
\end{thm}

\begin{proof}
Perform the replacements $a\to x$, $b\to-xa/A$, $c\to-q/a$, $d\to
AB/ab$, $e\to y$, $f\to -xa/B$ in \eqref{three-term} to achieve
  \bmn\label{three-term-a}
  &&\xxqdn{_8W_7}(x;-xa/A,-q/a,AB/ab,y,-xa/B;q,-bq/y)
  \nnm\\[1mm]\nnm
 &&\xxqdn\:=\frac{(xq,xabq/yAB,-bq/A,-Bq/ya,-ya,xa^2/B,-a/A,xya^2/AB;q)_{\infty}}
 {(xabq/AB,xq/y,-Bq/a,-bq/yA,-a,xya^2/B,-ya/A,xa^2/AB;q)_{\infty}}\\[1mm]
 &&\xxqdn\:\times\:{_8W_7}(xya^2/Bq;A,-xa^2b/AB,-ya/B,y,-xa/B;q,-B/b)
 \nnm\\[1mm]\nnm
  &&\xxqdn\:-\:\frac{(xq,-a/A,xa^2/A,-xa^2bq/A^2B,-xaq/yA,Bq/A,AB/ab,y,-xa/B;q)_{\infty}}
  {(-A/a,-xa,xabq/AB,xq/y,-Bq/a,-B/b,-ya/A,xa^2/AB,-yA/b;q)_{\infty}}\\[1mm]
  &&\xxqdn\:\times\:\frac{(A,ya/b,bq/ya;q)_{\infty}}{(-bq/yA,-a,xa^2q/A^2;q)_{\infty}}
  \nnm\\[1mm]
 &&\xxqdn\:\times\:
{_8W_7}(xa^2/A^2;-xa/A,q/A,-B/b,-ya/A,xa^2/AB;q,-bq/y).
 \emn
According to \eqref{two-term}, we have
 \bmn\label{two-term-a}
  &&\xxqdn{_8W_7}(xa^2/A^2;-xa/A,q/A,-B/b,-ya/A,xa^2/AB;q,-bq/y)
 \nnm\\[1mm]\nnm
 &&\xxqdn\:=
\frac{(xa^2q/A^2,-bq/A,-xa^2bq/yAB,Bq/y;q)_{\infty}}{(-xa^2bq/A^2B,Bq/A,xa^2q/yA,-bq/y;q)_{\infty}}\\
 &&\xxqdn\:\times\:{_8W_7}(xa^2/yA;q/y,-xa/y,xa^2/AB,-a,-B/b;q,-bq/A).
 \emn
Substitute \eqref{two-term-a} into \eqref{three-term-a} to attain

\bnm
  &&\xxqdn{_8W_7}(x;-xa/A,-q/a,AB/ab,y,-xa/B;q,-bq/y)
  \\[1mm]
 &&\xxqdn\:=\frac{(xq,xabq/yAB,-bq/A,-Bq/ya,-ya,xa^2/B,-a/A,xya^2/AB;q)_{\infty}}
 {(xabq/AB,xq/y,-Bq/a,-bq/yA,-a,xya^2/B,-ya/A,xa^2/AB;q)_{\infty}}\\[1mm]
 &&\xxqdn\:\times\:{_8W_7}(xya^2/Bq;A,-xa^2b/AB,-ya/B,y,-xa/B;q,-B/b)
 \\[1mm]
  &&\xxqdn\:-\:\frac{(xq,-xa/B,xa^2/A,-xaq/yA,-xa^2bq/yAB,A,-a/A,-bq/A;q)_{\infty}}
  {(-xa,xq/y,xabq/AB,xa^2/AB,xa^2q/yA,-a,-A/a,-Bq/a;q)_{\infty}}\\[1mm]
  &&\xxqdn\:\times\:\frac{(AB/ab,y,ya/b,bq/ya,Bq/y;q)_{\infty}}{(-B/b,-ya/A,-yA/b,-bq/yA,-bq/y;q)_{\infty}}
  \\[1mm]
 &&\xxqdn\:\times\:
{_8W_7}(xa^2/yA;q/y,-xa/y,xa^2/AB,-a,-B/b;q,-bq/A).
 \enm
Multiplying both sides by
$$\frac{1}{a}\frac{1-x}{(1+A/a)(1+B/a)(1-xab/AB)},$$
we obtain Theorem \ref{thm-a} after some simplification.
\end{proof}

\hspace{0.3cm}Secondly, we shall furnish a five-variable
generalization of \eqref{andrews} in the following theorem.

\begin{thm}  \label{thm-b}
Let $A$, $B$, $a$, $b$, $x$ be complex numbers with $|B|<|b|<1$.
Then
 \bnm
&&\xxqdn\xqdn\frac{1}{a}\sum_{k=0}^{\infty}\frac{(-q/a,AB/ab;q)_k(-b)^k}{(-A/a,-B/a;q)_{k+1}}
\frac{(x/A,x/B;q)_k(1+q^{2k}x/a)}{(x;q)_{k}(-xb/AB;q)_{k+1}}\\[1mm]
&&\xxqdn\xqdn\:-\frac{1}{b}\sum_{k=0}^{\infty}\frac{(A,-aq/B;q)_k(-B/b)^k}{(-a,-A/b;q)_{k+1}}
\frac{(xab/AB,x/B;q)_k(1+q^{2k}xa/B)}{(x;q)_{k}(-xa/AB;q)_{k+1}}\\[1mm]
&&\xxqdn\xqdn\:=\bigg(\frac{1}{a}-\frac{1}{b}\bigg)\frac{(q,aq/b,bq/a,A,B,AB/ab;q)_{\infty}}{(-a,-b,-A/a,-A/b,-B/a,-B/b;q)_{\infty}}\\
&&\xxqdn\xqdn\:\times\:\frac{(x/A,x/B,xab/AB;q)_{\infty}}{(x,-xa/AB,-xb/AB;q)_{\infty}}.
 \enm
\end{thm}

\begin{proof}
When $y=q$, it is routine to verify that
 \bnm
&&\lambda(y)=1,\\[1mm]
&&\theta(y)=-\frac{1}{A+b},\\[1mm]
&&{_8W_7}(xa^2/yA;q/y,-xa/y,xa^2/AB,-a,-B/b;q,-bq/A)=1.
 \enm
So the case $y=q$ of Theorem \ref{thm-a} reads
 \bnm
&&\xqdn\frac{1}{a}\sum_{k=0}^{\infty}\frac{(-q/a,AB/ab;q)_k(-b)^k}{(-A/a,-B/a;q)_{k+1}}
\frac{(-xa/A,-xa/B;q)_k(1-q^{2k}x)}{(-xa;q)_{k}(xab/AB;q)_{k+1}}\\[1mm]
&&\xqdn\:-\frac{1}{b}\sum_{k=0}^{\infty}\frac{(A,-aq/B;q)_k(-B/b)^k}{(-a,-A/b;q)_{k+1}}
\frac{(-xa^2b/AB,-xa/B;q)_k(1-q^{2k}xa^2/B)}{(-xa;q)_{k}(xa^2/AB;q)_{k+1}}\\[1mm]
&&\xqdn\:=\bigg(\frac{1}{a}-\frac{1}{b}\bigg)\frac{(q,aq/b,bq/a,A,B,AB/ab;q)_{\infty}}{(-a,-b,-A/a,-A/b,-B/a,-B/b;q)_{\infty}}\\[1mm]
&&\xqdn\:\times\:\frac{(-xa/A,-xa/B,-xa^2b/AB;q)_{\infty}}{(-xa,xa^2/AB,xab/AB;q)_{\infty}}.
 \enm
Replacing $x$ by $-x/a$ in the last equation, we get Theorem
\ref{thm-b}.
\end{proof}

\hspace{0.3cm}When $x=0$, Theorem \ref{thm-b} reduces to Andrews'
identity \eqref{andrews} exactly. Hence Theorem \ref{thm-b} is a
five-variable generalization of \eqref{ramanujan} and Theorem
\ref{thm-a} can be regarded as a six-variable generalization of
\eqref{ramanujan}.

\section{A six-variable generalization  of Ramanujan's\\ reciprocity theorem}

\hspace{0.3cm} Firstly, we shall establish a six-variable
generalization of Ramanujan's reciprocity theorem
\eqref{Ramanujan-a}, which is different from \eqref{liu-b}, in the
following theorem.

\begin{thm}  \label{thm-c}
Let $a$, $c$, $u$, $v$, $x$, $y$ be complex numbers. Then
 \bnm
&&\xqdn
v\sum_{k=0}^{\infty}\frac{(q/u,acuv;q)_kv^k}{(av,cv;q)_{k+1}}
\frac{(x,y,x/av,x/cv;q)_k(1-q^{2k}x)}{(q,xu,xq/y;q)_{k}(x/acuv;q)_{k+1}}\bigg(\frac{q}{y}\bigg)^k\\[1mm]
&&\xqdn\:-\:u\rho(y)\sum_{k=0}^{\infty}\frac{(y/v,acuv;q)_ku^k}{(yau/q,ycu/q;q)_{k+1}}
\frac{(xyu/vq,y,x/av,x/cv;q)_k(1-q^{2k-1}xyu/v)}{(q,xu,xu/v;q)_{k}(xy/acv^2q;q)_{k+1}}\\[1mm]
&&\xqdn\:=\frac{\omega(y)(acuv,cuv,y,yu/v,vq/yu,acuvq/y,x/av,xq/ycv,xq/yacv;q)_{\infty}}
{(u,au,av,cv,xu,ycu,vq/y,x/acuv,x/acv^2;q)_{\infty}}\\[1mm]
&&\xqdn\:\times\:{_8W_7}(xu/ycv;q/y,xu/y,x/acv^2,u,au;q,q/cu),
 \enm
where $\max\{|qv/y|,|u|,|q/cu|\}<1$ and the symbols $\rho(y)$,
$\omega(y)$ designate
 \bnm
  &&\xqdn\rho(y)=\frac{y}{q}\frac{(x,xu/av,q/cu,q/cv,yu/q,avq/y,xq/yacuv,xy/acv^2q;q)_{\infty}}
  {(xq/y,xyu/avq,q^2/ycu,y/cv,u,av,x/acuv,x/acv^2;q)_{\infty}},\\[1mm]
   &&\xqdn\omega(y)=\frac{u}{cu-1}\frac{(x,xu/cv,1/cu,q/cv;q)_{\infty}}
  {(xq/y,xuq/ycv,q/ycu,y/cv;q)_{\infty}}.
  \enm
\end{thm}

\begin{proof}
Let $\Omega(A,B,a,b,x,y)$ denote the second series on the left hand
side of Theorem \ref{thm-a}. Then we have
 \bmn\label{relation-a}
\xqdn\Omega(A,B,a,b,x,y)&&\xqdn\!=\sum_{k=0}^{\infty}\frac{(A,-ya/B;q)_k(-B/b)^k}{(-ya/q,-yA/bq;q)_{k+1}}
\frac{(xya^2/Bq,y,-xa^2b/AB;q)_k}{(q,-xa,xa^2/B;q)_{k}}
 \nnm\\[1mm]\nnm
&&\!\xqdn\times\:\frac{(-xa/B;q)_k(1-q^{2k-1}xya^2/B)}{(xya^2/ABq;q)_{k+1}}
\nnm\\[1mm]\nnm
&&\!\xqdn=\frac{1-xya^2/Bq}{(1+ya/q)(1+yA/bq)(1-xya^2/ABq)}
 \\[1mm]
&&\!\xqdn\times\:{_8W_7}(xya^2/Bq;A,-xa^2b/AB,-ya/B,y,-xa/B;q,-B/b).
 \emn
In terms of \eqref{two-term}, we gain the relation:
  \bmn\label{relation-b}
  &&\xxqdn{_8W_7}(xya^2/Bq;A,-xa^2b/AB,-ya/B,y,-xa/B;q,-B/b)
  \nnm\\[1mm]\nnm
 &&\xxqdn\:=
\frac{(xya^2/B,-a,xa/b,-yB/b;q)_{\infty}}{(xa^2/B,-ya,xya/b,-B/b;q)_{\infty}}\nnm\\[1mm]
&&\xxqdn\:\times\:{_8W_7}(xya/bq;AB/ab,-xa/A,-y/b,y,-xa/B;q,-a).
 \emn
The combination of \eqref{relation-a} and \eqref{relation-b}
produces
 \bnm
\Omega(A,B,a,b,x,y)&&\xqdn\!=\frac{(xya^2/Bq,-a,xa/b,-yB/bq;q)_{\infty}}{(xa^2/B,-ya/q,xya/bq,-B/b;q)_{\infty}}
 \enm
 \bmn\label{relation-c}
&&\!\xqdn\times\:\sum_{k=0}^{\infty}\frac{(-y/b,AB/ab;q)_k(-a)^k}{(-yA/bq,-yB/bq;q)_{k+1}}
\frac{(xya/bq,y,-xa/A,-xa/B;q)_k}{(q,-xa,xa/b;q)_{k}}
\nnm\\[1mm]
&&\xqdn\:\times\:\frac{(1-q^{2k-1}xya/b)}{(xya^2/ABq;q)_{k+1}}.
 \emn
Substitute \eqref{relation-c} into Theorem \ref{thm-a} to achieve
 \bnm
&&\xqdn\frac{1}{a}\sum_{k=0}^{\infty}\frac{(-q/a,AB/ab;q)_k(-b)^k}{(-A/a,-B/a;q)_{k+1}}
\frac{(x,y,-xa/A,-xa/B;q)_k(1-q^{2k}x)}{(q,-xa,xq/y;q)_{k}(xab/AB;q)_{k+1}}\bigg(\frac{q}{y}\bigg)^k\\[1mm]
&&\xqdn\:-\frac{\eta(y)}{b}\sum_{k=0}^{\infty}\frac{(-y/b,AB/ab;q)_k(-a)^k}{(-yA/bq,-yB/bq;q)_{k+1}}
\frac{(xya/bq,y,-xa/A,-xa/B;q)_k}{(q,-xa,xa/b;q)_{k}}\\[1mm]
&&\xqdn\:\times\:\frac{(1-q^{2k-1}xya/b)}{(xya^2/ABq;q)_{k+1}}\\[1mm]
&&\xqdn\:=\frac{\theta(y)(y,ya/b,bq/ya,A,Bq/y,AB/ab,-xa/B,-xaq/yA,-xa^2bq/yAB;q)_{\infty}}{(-a,-bq/y,-A/a,-yA/b,-B/a,-B/b,-xa,xab/AB,xa^2/AB;q)_{\infty}}\\[1mm]
&&\xqdn\:\times\:{_8W_7}(xa^2/yA;q/y,-xa/y,xa^2/AB,-a,-B/b;q,-bq/A),
 \enm
where $\theta(y)$ has appeared in Theorem \ref{thm-a} and the
notation $\eta(y)$ stands for
 \bnm
  &&\xqdn\eta(y)=\frac{y}{q}\frac{(x,xa/b,-aq/A,-bq/A,-Bq/ya,-yB/bq,xabq/yAB,xya^2/ABq;q)_{\infty}}
  {(xq/y,xya/bq,-ya/A,-bq^2/yA,-B/a,-B/b,xab/AB,xa^2/AB;q)_{\infty}}.
  \enm
Employing the replacements $A\to cuv$, $B\to auv$, $a\to-u$, $b\to
-v$ in the last equation, we attain Theorem \ref{thm-c} after some
simplification.
\end{proof}

\hspace{0.3cm}Secondly, we shall offer a five-variable
generalization of \eqref{Ramanujan-a} in the following theorem.

\begin{thm}  \label{thm-d}
Let $a$, $c$, $u$, $v$, $x$ be complex numbers with
$\max\{|u|,|v|\}<1$. Then
 \bnm
&&\xqdn
v\sum_{k=0}^{\infty}\frac{(q/u,acuv;q)_kv^k}{(av,cv;q)_{k+1}}
\frac{(x/a,x/c;q)_k(1-q^{2k}xv)}{(xuv;q)_{k}(x/acu;q)_{k+1}}\\[1mm]
&&\xqdn\:-\:u\sum_{k=0}^{\infty}\frac{(q/v,acuv;q)_ku^k}{(au,cu;q)_{k+1}}
\frac{(x/a,x/c;q)_k(1-q^{2k}xu)}{(xuv;q)_{k}(x/acv;q)_{k+1}}\\[1mm]
&&\xqdn\:=\frac{(v-u)(q,uq/v,vq/u,auv,cuv,acuv;q)_{\infty}}{(u,v,au,av,cu,cv;q)_{\infty}}
\frac{(x/a,x/c,x/ac;q)_{\infty}}{(xuv,x/acu,x/acv;q)_{\infty}}.
  \enm
\end{thm}

\begin{proof}
When $y=q$, it is easy to see that
 \bnm
&&\rho(y)=1,\\[1mm]
&&\omega(y)=\frac{u}{cu-1},\\[1mm]
&&{_8W_7}(xu/ycv;q/y,xu/y,x/acv^2,u,au;q,q/cu)=1.
 \enm
Therefore, the case $y=q$ of Theorem \ref{thm-c} reads

 \bnm
&&\xqdn
v\sum_{k=0}^{\infty}\frac{(q/u,acuv;q)_kv^k}{(av,cv;q)_{k+1}}
\frac{(x/av,x/cv;q)_k(1-q^{2k}x)}{(xu;q)_{k}(x/acuv;q)_{k+1}}\\[1mm]
&&\xqdn\:-\:u\sum_{k=0}^{\infty}\frac{(q/v,acuv;q)_ku^k}{(au,cu;q)_{k+1}}
\frac{(x/av,x/cv;q)_k(1-q^{2k}xu/v)}{(xu;q)_{k}(x/acv^2;q)_{k+1}}\\[1mm]
&&\xqdn\:=\frac{(v-u)(q,uq/v,vq/u,auv,cuv,acuv;q)_{\infty}}{(u,v,au,av,cu,cv;q)_{\infty}}
\frac{(x/av,x/cv,x/acv;q)_{\infty}}{(xu,x/acuv,x/acv^2;q)_{\infty}}.
  \enm
Replacing $x$ by $xv$ in the last equation, we obtain Theorem
\ref{thm-d}.
\end{proof}

\hspace{0.3cm}When $x=0$, Theorem \ref{thm-d} reduces to
\eqref{liu-a} exactly. Thus Theorem \ref{thm-d} is a five-variable
generalization of \eqref{Ramanujan-a} and Theorem \ref{thm-c} can be
regarded as a six-variable generalization of  \eqref{Ramanujan-a}.

\section{Some transformation formulas involving bilateral basic
hypergeometric series}
\hspace{0.3cm}Above all, we shall give two limiting relations in the
following two lemmas.

\begin{lemm}  \label{lemma-a}
Let $a$ and $\{b_i\}_{i=1}^{2r+1}$ be complex numbers. Then
 \bnm
&&\lim_{a\to1}{_{2r+4}W_{2r+3}}\Big(a;b_1,b_2,\ldots,b_{2r+1};q,a^rq^r/b_1b_2\cdots
b_{2r+1}\Big)\\[1mm]
&&\:={_{2r+1}\psi_{2r+1}}\ffnk{cccccc}{q,\frac{q^r}{b_1b_2\cdots
b_{2r+1}}}{b_1,b_2,\ldots,b_{2r+1}}{q/b_1,q/b_2,\ldots,q/b_{2r+1}},
  \enm
provided $\max\{|a^rq^r/b_1b_2\cdots b_{2r+1}|,|q^r/b_1b_2\cdots
b_{2r+1}|\}<1$.
\end{lemm}

\begin{proof}
 \bnm
&&\lim_{a\to1}{_{2r+4}W_{2r+3}}\Big(a;b_1,b_2,\ldots,b_{2r+1};q,a^rq^r/b_1b_2\cdots
b_{2r+1}\Big)\\[1mm]
&&=
 1+\:\sum_{k=1}^{\infty}\{1+q^k\}\frac{(b_1,b_2,\ldots,b_{2r+1};q)_k}{(q/b_1,q/b_2,\ldots,q/b_{2r+1};q)_k}
 \bigg(\frac{q^r}{b_1b_2\cdots b_{2r+1}}\bigg)^k\\
&&=1+\:\sum_{k=1}^{\infty}\frac{(b_1,b_2,\ldots,b_{2r+1};q)_k}{(q/b_1,q/b_2,\ldots,q/b_{2r+1};q)_k}
 \bigg(\frac{q^r}{b_1b_2\cdots b_{2r+1}}\bigg)^k\\
&&\quad\:\:\:+\:\sum_{k=1}^{\infty}
 \frac{(b_1,b_2,\ldots,b_{2r+1};q)_k}{(q/b_1,q/b_2,\ldots,q/b_{2r+1};q)_k}
 \bigg(\frac{q^{r+1}}{b_1b_2\cdots b_{2r+1}}\bigg)^k\\
&&=1+\:\sum_{k=1}^{\infty}\frac{(b_1,b_2,\ldots,b_{2r+1};q)_k}{(q/b_1,q/b_2,\ldots,q/b_{2r+1};q)_k}
 \bigg(\frac{q^r}{b_1b_2\cdots b_{2r+1}}\bigg)^k\\
&&\quad\:\:\:+\sum_{k=-\infty}^{-1}
 \frac{(b_1,b_2,\ldots,b_{2r+1};q)_{-k}}{(q/b_1,q/b_2,\ldots,q/b_{2r+1};q)_{-k}}
 \bigg(\frac{q^{r+1}}{b_1b_2\cdots b_{2r+1}}\bigg)^{-k}\\
&&=1+\:\sum_{k=1}^{\infty}\frac{(b_1,b_2,\ldots,b_{2r+1};q)_k}{(q/b_1,q/b_2,\ldots,q/b_{2r+1};q)_k}
 \bigg(\frac{q^r}{b_1b_2\cdots b_{2r+1}}\bigg)^k
 \enm

 \bnm
&&\quad\:\:\:+\sum_{k=-\infty}^{-1}
 \frac{(b_1,b_2,\ldots,b_{2r+1};q)_k}{(q/b_1,q/b_2,\ldots,q/b_{2r+1};q)_k}
 \bigg(\frac{q^r}{b_1b_2\cdots b_{2r+1}}\bigg)^k\\
 \\&&=\:{_{2r+1}\psi_{2r+1}}\ffnk{cccccc}{q,\frac{q^r}{b_1b_2\cdots
b_{2r+1}}}{b_1,b_2,\ldots,b_{2r+1}}{q/b_1,q/b_2,\ldots,q/b_{2r+1}}.
 \enm
This completes the proof of Lemma \ref{lemma-a}.
\end{proof}

\begin{lemm}  \label{lemma-b}
Let $a$ and $\{b_i\}_{i=1}^{2r+1}$ be complex numbers. Then
 \bnm
&&\lim_{a\to{q}}{_{2r+4}W_{2r+3}}\Big(a;b_1,b_2,\ldots,b_{2r+1};q,a^rq^r/b_1b_2\cdots
b_{2r+1}\Big)\\[1mm]
&&\:=\frac{1}{1-q}\,{_{2r+1}\psi_{2r+1}}\ffnk{cccccc}{q,\frac{q^{2r}}{b_1b_2\cdots
b_{2r+1}}}{b_1,b_2,\ldots,b_{2r+1}}{q^2/b_1,q^2/b_2,\ldots,q^2/b_{2r+1}},
  \enm
where $\max\{|a^rq^r/b_1b_2\cdots b_{2r+1}|,|q^{2r}/b_1b_2\cdots
b_{2r+1}|\}<1$.
\end{lemm}

\begin{proof}
 \bnm
&&\lim_{a\to{q}}{_{2r+4}W_{2r+3}}\Big(a;b_1,b_2,\ldots,b_{2r+1};q,a^rq^r/b_1b_2\cdots
b_{2r+1}\Big)\\[1mm]
&&=
 \frac{1}{1-q}\sum_{k=0}^{\infty}\{1-q^{1+2k}\}\frac{(b_1,b_2,\ldots,b_{2r+1};q)_k}{(q^2/b_1,q^2/b_2,\ldots,q^2/b_{2r+1};q)_k}
 \bigg(\frac{q^{2r}}{b_1b_2\cdots b_{2r+1}}\bigg)^k\\
&&=\frac{1}{1-q}\sum_{k=0}^{\infty}\frac{(b_1,b_2,\ldots,b_{2r+1};q)_k}{(q^2/b_1,q^2/b_2,\ldots,q^2/b_{2r+1};q)_k}
 \bigg(\frac{q^{2r}}{b_1b_2\cdots b_{2r+1}}\bigg)^k\\
&&-\:\frac{1}{1-q}\sum_{k=0}^{\infty}\frac{(b_1,b_2,\ldots,b_{2r+1};q)_k}{(q^2/b_1,q^2/b_2,\ldots,q^2/b_{2r+1};q)_k}
 \frac{q^{2(r+1)k+1}}{(b_1b_2\cdots b_{2r+1})^k}\\
&&=\frac{1}{1-q}\sum_{k=0}^{\infty}\frac{(b_1,b_2,\ldots,b_{2r+1};q)_k}{(q^2/b_1,q^2/b_2,\ldots,q^2/b_{2r+1};q)_k}
 \bigg(\frac{q^{2r}}{b_1b_2\cdots b_{2r+1}}\bigg)^k\\
&&-\:\frac{1}{1-q}\sum_{k=-\infty}^{-1}\frac{(b_1,b_2,\ldots,b_{2r+1};q)_{-k-1}}{(q^2/b_1,q^2/b_2,\ldots,q^2/b_{2r+1};q)_{-k-1}}
 \frac{q^{2(r+1)(-k-1)+1}}{(b_1b_2\cdots b_{2r+1})^{-k-1}}\\
&&=\frac{1}{1-q}\sum_{k=0}^{\infty}\frac{(b_1,b_2,\ldots,b_{2r+1};q)_k}{(q^2/b_1,q^2/b_2,\ldots,q^2/b_{2r+1};q)_k}
 \bigg(\frac{q^{2r}}{b_1b_2\cdots b_{2r+1}}\bigg)^k\\
&&+\:\frac{1}{1-q}\sum_{k=-\infty}^{-1}\frac{(b_1,b_2,\ldots,b_{2r+1};q)_k}{(q^2/b_1,q^2/b_2,\ldots,q^2/b_{2r+1};q)_k}
 \bigg(\frac{q^{2r}}{b_1b_2\cdots b_{2r+1}}\bigg)^k\\
&&\:=\frac{1}{1-q}\,{_{2r+1}\psi_{2r+1}}\ffnk{cccccc}{q,\frac{q^{2r}}{b_1b_2\cdots
b_{2r+1}}}{b_1,b_2,\ldots,b_{2r+1}}{q^2/b_1,q^2/b_2,\ldots,q^2/b_{2r+1}}.
 \enm
This finishes the proof of Lemma \ref{lemma-b}.
\end{proof}

\hspace{0.3cm}Let the symbol ``idem$(x,y)$'' after an expression
mean that the preceding expression is repeated with $x$ and $y$
interchanged. Subsequently, some transformation formulas involving
bilateral basic hypergeometric series will be furnished in the
following theorems.

\begin{thm}  \label{thm-e}
Let $b,c,d,e,f$ be complex numbers. Then
 \bmn\label{bailey-gen-a}
&&\xxqdn{_5\psi_5}\ffnk{cccccc}{q,\frac{q^{2}}{bcdef}}{b,c,d,e,f}{q/b,q/c,q/d,q/e,q/f}
\nnm\\[1mm]\nnm
&&\xxqdn\:=\frac{(q,q/be,q/ce,q/de,eq/c,eq/d,f,f,q^2/cdf,q^2/bcdf;q)_{\infty}}
{(q/b,q/c,q/c,q/d,q/d,q/e,f/e,ef,eq^2/cdf,q^2/bcdef;q)_{\infty}}\\[1mm]
&&\xxqdn\:\times\:{_8W_7}(eq/cdf;e,be,q/cd,q/cf,q/df;q,q/b)+idem(e;f),
  \emn
provided $\max\{|q/b|,|q^2/bcdef|\}<1$,
 \bmn\label{bailey-gen-b}
&&{_5\psi_5}\ffnk{cccccc}{q,\frac{q^{4}}{bcdef}}{b,c,d,e,f}{q^2/b,q^2/c,q^2/d,q^2/e,q^2/f}
\nnm\\[1mm]\nnm
&&\:=\frac{(q,q^2/be,q^2/ce,q^2/de,eq/c,eq/d,f,f/q,q^3/cdf,q^4/bcdf;q)_{\infty}}
{(q^2/b,q/c,q^2/c,q/d,q^2/d,q^2/e,f/e,ef/q,eq^3/cdf,q^4/bcdef;q)_{\infty}}\\[1mm]
&&\:\times\:{_8W_7}(eq^2/cdf;e,be/q,q^2/cd,q^2/cf,q^2/df;q,q/b)+idem(e;f),
  \emn
where $\max\{|q/b|,|q^4/bcdef|\}<1$.
\end{thm}

\begin{proof}
Performing the replacements $b\to f/a, c\to e/a, d\to c/a, e\to d/a,
f\to b/a$ in \eqref{bailey-three-term} and then replacing $a$ by
$a^{1/2}$, we get
 \bmn\label{equation-c}
  &&\xqdn{_8W_7}(a;b,c,d,e,f;q,a^2q^2/bcdef)
  \nnm\\[1mm]\nnm
  &&\xqdn\:=
\frac{(aq,aq/be,aq/ce,aq/de,eq/a,eq/b,eq/c,eq/d,f,f/a;q)_{\infty}}
{(q/b,q/c,q/d,aq/b,aq/c,aq/d,aq/e,f/e,ef/a,e^2q/a;q)_{\infty}}
 \\[1mm]&&\xqdn\:\times\:
 {_8W_7}(e^2/a;e,be/a,ce/a,de/a,ef/a;q,a^2q^2/bcdef)+idem(e,f).
 \emn
By means of \eqref{two-term}, we have
 \bnm
  &&\xqdn{_8W_7}(e^2/a;e,be/a,ce/a,de/a,ef/a;q,a^2q^2/bcdef)
  \nnm\\[1mm]\nnm
  &&\:\xqdn\:=
\frac{(e^2q/a,q/b,aq^2/cdf,a^2q^2/bcdf;q)_{\infty}}
{(eq/a,eq/b,aeq^2/cdf,a^2q^2/bcdef;q)_{\infty}}
 \\[1mm]&&\xqdn\,\:\times\:
 {_8W_7}(aeq/cdf;e,be/a,aq/cd,aq/cf,aq/df;q,q/b),
 \\[1mm]
  &&\xqdn{_8W_7}(f^2/a;f,bf/a,cf/a,df/a,ef/a;q,a^2q^2/bcdef)
  \nnm\\[1mm]\nnm
  &&\:\xqdn\:=
\frac{(f^2q/a,q/b,aq^2/cde,a^2q^2/bcde;q)_{\infty}}
{(fq/a,fq/b,afq^2/cde,a^2q^2/bcdef;q)_{\infty}}
 \\[1mm]&&\xqdn\,\:\times\:
 {_8W_7}(afq/cde;f,bf/a,aq/cd,aq/ce,aq/de;q,q/b).
 \enm
Substitute the last two relations into \eqref{equation-c} to derive
 \bmn\label{equation-d}
  &&\xqdn{_8W_7}(a;b,c,d,e,f;q,a^2q^2/bcdef)
  \nnm\\[1mm]\nnm
  &&\xqdn\:=
\frac{(aq,aq/be,aq/ce,aq/de,eq/c,eq/d,f,f/a,aq^2/cdf,a^2q^2/bcdf;q)_{\infty}}
{(q/c,q/d,aq/b,aq/c,aq/d,aq/e,f/e,ef/a,aeq^2/cdf,a^2q^2/bcdef;q)_{\infty}}
 \\[1mm]&&\xqdn\:\times\:
 {_8W_7}(aeq/cdf;e,be/a,aq/cd,aq/cf,aq/df;q,q/b)+idem(e,f).
 \emn
Letting $a\to 1$ in \eqref{equation-d} and using Lemma
\ref{lemma-a}, we gain \eqref{bailey-gen-a}. Letting $a\to q$ in
\eqref{equation-d} and utilizing Lemma \ref{lemma-b}, we achieve
\eqref{bailey-gen-b}.
\end{proof}

\hspace{0.3cm}It is obvious that \eqref{bailey-gen-a} is a
generalization of \eqref{p33-a}. When $df=q$, The former reduces to
the latter. It is also apparent that \eqref{bailey-gen-b} is a
generalization of \eqref{p33-b}. When $df=q^2$, The former reduces
to the latter.

\begin{thm}  \label{thm-f}
Let $a,b,c,d$ be complex numbers. Then
 \bnm
&&\xqdn\xxqdn{_5\psi_5}\ffnk{cccccc}{q,-q}{a^{\frac{1}{2}},-a^{\frac{1}{2}},b,c,d}{qa^{-\frac{1}{2}},-qa^{-\frac{1}{2}},q/b,q/c,q/d}\\[1mm]
&&\xqdn\xxqdn\:=\frac{(q,-q/a,qa^{\frac{1}{2}},-qa^{\frac{1}{2}};q)_{\infty}}
{(aq,-q,qa^{-\frac{1}{2}},-qa^{-\frac{1}{2}};q)_{\infty}}
{_4\phi_3}\ffnk{cccccc}{q,-\frac{q}{a}}{a,ab,ac,ad}{q/b,q/c,q/d},
  \enm
provided $q=abcd$ and $|q/a|<1$,
  \bnm
&&{_5\psi_5}\ffnk{cccccc}{q,-q}{a^{\frac{1}{2}},-a^{\frac{1}{2}},b,c,d}{q^2a^{-\frac{1}{2}},-q^2a^{-\frac{1}{2}},q^2/b,q^2/c,q^2/d}\\[1mm]
&&\:=\frac{(q,-q^2/a,qa^{\frac{1}{2}},-qa^{\frac{1}{2}};q)_{\infty}}
{(aq,-q,q^2a^{-\frac{1}{2}},-q^2a^{-\frac{1}{2}};q)_{\infty}}
{_4\phi_3}\ffnk{cccccc}{q,-\frac{q^2}{a}}{a,ab/q,ac/q,ad/q}{q^2/b,q^2/c,q^2/d},
  \enm
where $q^3=abcd$ and $|q^2/a|<1$.
\end{thm}

\begin{proof}
Employ the replacements $b\to ab, c\to ac, d\to ad$ in
\eqref{equation-a} to attain
 \bnm
  &&\xxqdn{_8W_7}\Big(\mu;a^{\frac{1}{2}},-a^{\frac{1}{2}},\mu b,\mu c,\mu d;q,-q\Big)
  \nnm\\[1mm]
  &&\xxqdn\:=\frac{(\mu q,-\mu q/a,qa^{\frac{1}{2}},-qa^{\frac{1}{2}};q)_{\infty}}{(aq,-q,\mu qa^{-\frac{1}{2}},-\mu qa^{-\frac{1}{2}};q)_{\infty}}
  {_4\phi_3}\ffnk{ccccccc}{q,-\frac{\mu q}{a}}{a,ab,ac,ad}{q/b,q/c,q/d},
  \enm
provided $\mu=q/abcd$. Letting $\mu\to 1$ in the last equation and
using Lemma \ref{lemma-a}, we obtain the first equation of Theorem
\ref{thm-f}.

\hspace{0.3cm}Perform the replacements $b\to ab/q, c\to ac/q, d\to
ad/q$ in \eqref{equation-a} to get
 \bnm
  &&\qdn{_8W_7}\Big(\nu;a^{\frac{1}{2}},-a^{\frac{1}{2}},\nu b/q,\nu c/q,\nu d/q;q,-q\Big)
  \nnm\\[1mm]
  &&\qdn\:=\frac{(\nu q,-\nu q/a,qa^{\frac{1}{2}},-qa^{\frac{1}{2}};q)_{\infty}}{(aq,-q,\nu qa^{-\frac{1}{2}},-\nu qa^{-\frac{1}{2}};q)_{\infty}}
  {_4\phi_3}\ffnk{ccccccc}{q,-\frac{\nu q}{a}}{a,ab/q,ac/q,ad/q}{q^2/b,q^2/c,q^2/d},
  \enm
where $\nu=q^4/abcd$. Letting $\nu\to q$ in the last equation and
utilizing Lemma \ref{lemma-b}, we gain the second equation of
Theorem \ref{thm-f}.
\end{proof}

\begin{thm}  \label{thm-g}
Let $a,b,c,d$ be complex numbers. Then
 \bnm
&&\xqdn\xxqdn{_7\psi_7}\ffnk{cccccc}{q,\frac{q}{a}}{a^{\frac{1}{2}},-a^{\frac{1}{2}},(qa)^{\frac{1}{2}},-(qa)^{\frac{1}{2}},b,c,d}
{qa^{-\frac{1}{2}},-qa^{-\frac{1}{2}},(q/a)^{\frac{1}{2}},-(q/a)^{\frac{1}{2}},q/b,q/c,q/d}\\[1mm]
&&\xqdn\xxqdn\:=\frac{(q,q/a^2;q)_{\infty}}{(q/a,q/a;q)_{\infty}}
{_4\phi_3}\ffnk{cccccc}{q,\frac{q}{a^2}}{a,ab,ac,ad}{q/b,q/c,q/d},
  \enm
provided $q=abcd$ and $\max\{|q/a|,|q/a^2|\}<1$,
 \bnm
&&\qqdn{_7\psi_7}\ffnk{cccccc}{q,\frac{q^2}{a}}{a^{\frac{1}{2}},-a^{\frac{1}{2}},(qa)^{\frac{1}{2}},-(qa)^{\frac{1}{2}},b,c,d}
{q^2a^{-\frac{1}{2}},-q^2a^{-\frac{1}{2}},(q^3/a)^{\frac{1}{2}},-(q^3/a)^{\frac{1}{2}},q^2/b,q^2/c,q^2/d}\\[1mm]
&&\qqdn\:=\frac{(q,q^3/a^2;q)_{\infty}}{(q^2/a,q^3/a;q)_{\infty}}
{_4\phi_3}\ffnk{cccccc}{q,\frac{q^3}{a^2}}{a,ab/q,ac/q,ad/q}{q^2/b,q^2/c,q^2/d},
  \enm
where $q^3=abcd$ and $\max\{|q^2/a|,|q^3/a^2|\}<1$.
\end{thm}

\begin{proof}
Employ the replacements $b\to ab, c\to ac, d\to ad$ in
\eqref{equation-b} to achieve
 \bnm
  &&\xxqdn{_{10}W_9}\Big(\mu;a^{\frac{1}{2}},-a^{\frac{1}{2}},(aq)^{\frac{1}{2}},-(aq)^{\frac{1}{2}},\mu b,\mu c,\mu d;q,\mu q/a\Big)
  \nnm\\[1mm]
  &&\xxqdn\:=\frac{(\mu q,\mu^2 q/a^2;q)_{\infty}}{(\mu q/a, \mu^2 q/a;q)_{\infty}}
  {_4\phi_3}\ffnk{ccccccc}{q,\frac{\mu^2q}{a^2}}{a,ab,ac,ad}{q/b,q/c,q/d}.
  \enm
 Letting $\mu\to 1$ in the last equation and using Lemma
\ref{lemma-a}, we attain the first equation of Theorem \ref{thm-g}.

\hspace{0.3cm}Perform the substitutions $b\to ab/q, c\to ac/q, d\to
ad/q$ in \eqref{equation-b} to obtain
  \bnm
  &&{_{10}W_9}\Big(\nu;a^{\frac{1}{2}},-a^{\frac{1}{2}},(aq)^{\frac{1}{2}},-(aq)^{\frac{1}{2}},\nu b/q,\nu c/q,\nu d/q;q,\nu q/a\Big)
  \nnm\\[1mm]
  &&\:=\frac{(\nu q,\nu^2 q/a^2;q)_{\infty}}{(\nu q/a, \nu^2 q/a;q)_{\infty}}
  {_4\phi_3}\ffnk{ccccccc}{q,\frac{\mu^2q}{a^2}}{a,ab/q,ac/q,ad/q}{q^2/b,q^2/c,q^2/d}.
  \enm
Letting $\nu\to q$ in the last equation and utilizing Lemma
\ref{lemma-b}, we get the second equation of Theorem \ref{thm-g}.
\end{proof}

By applying Lemmas \ref{lemma-a} and \ref{lemma-b} to
\eqref{Guo-Schlosser}, we can deduce the following results.

\begin{thm}  \label{thm-h}
Let $b,c,d$ be complex numbers . Then
 \bnm
&&\qqdn\xqdn\xxqdn{_9\psi_9}\ffnk{cccccc}{q,\frac{q}{b}}{(q/b)^{\frac{1}{2}},-(q/b)^{\frac{1}{2}},(q^2/b)^{\frac{1}{2}},-(q^2/b)^{\frac{1}{2}},b,c,d,b/c,b/d}
{(qb)^{\frac{1}{2}},-(qb)^{\frac{1}{2}},b^{\frac{1}{2}},-b^{\frac{1}{2}},q/b,q/c,q/d,cq/b,dq/b}\\[1mm]
&&\qqdn\xqdn\xxqdn\:=\frac{(q,b/c,b/d,q/cd;q)_{\infty}}{(b,q/c,q/d,b/cd;q)_{\infty}}
{_4\phi_3}\ffnk{cccccc}{q,\frac{q^2}{b^2}}{b,c,d,cd}{cq/b,dq/b,cdq/b}\\
&&\qqdn\xqdn\xxqdn\:+\,\frac{(q,c,d,cdq/b^2;q)_{\infty}}{(b,cq/b,dq/b,cd/b;q)_{\infty}}
{_4\phi_3}\ffnk{cccccc}{q,\frac{q^2}{b^2}}{b,b/c,b/d,b^2/cd}{q/c,q/d,bq/cd},
  \enm
provided $|q/b|<1$,
 \bnm
&&\xqdn{_9\psi_9}\ffnk{cccccc}{q,\frac{q}{b}}{(q^2/b)^{\frac{1}{2}},-(q^2/b)^{\frac{1}{2}},(q^3/b)^{\frac{1}{2}},-(q^3/b)^{\frac{1}{2}},b,c,d,bq/c,bq/d}
{(q^2b)^{\frac{1}{2}},-(q^2b)^{\frac{1}{2}},(qb)^{\frac{1}{2}},-(qb)^{\frac{1}{2}},q^2/b,q^2/c,q^2/d,cq/b,dq/b}\\[1mm]
&&\xqdn\:=\frac{(q,bq/c,bq/d,q^2/cd;q)_{\infty}}{(bq,q^2/c,q^2/d,bq/cd;q)_{\infty}}
{_4\phi_3}\ffnk{cccccc}{q,\frac{q^2}{b^2}}{b,c,d,cd/q}{cq/b,dq/b,cd/b}\\
&&\xqdn\:+\,\frac{(q,c,d,cd/b^2;q)_{\infty}}{(bq,cq/b,dq/b,cd/bq;q)_{\infty}}
{_4\phi_3}\ffnk{cccccc}{q,\frac{q^2}{b^2}}{b,bq/c,bq/d,b^2q/cd}{q^2/c,q^2/d,bq^2/cd},
\enm
 where $|q/b|<1$.
\end{thm}

\textbf{Acknowledgments}

 The work is supported by the National Natural Science Foundations of China (Nos. 11661032,11601151).



\end{document}